\documentclass[11pt]{article}
\usepackage{amsmath}
\usepackage{amssymb}
\usepackage{eucal}
\usepackage{stmaryrd}

\begin{document}

\baselineskip 16pt

\title{Finite groups with  systems of  $K$-$\frak{F}$-subnormal   subgroups }

\author            
{Vladimir N. Semenchuk, Alexander  N. Skiba \\
{\small Department of Mathematics and Technologies of Programming, 
 Francisk Skorina Gomel State University,}\\
{\small Gomel 246019, Belarus}\\
{\small E-mail:   semenchuk@gsu.by, alexander.skiba49@gmail.com}}

\date{}
\maketitle

\begin{abstract} Let  $\frak {F}$  be  a  class of group. A 
  subgroup $A$ of  a finite group $G$
 is said to be  \emph{$K$-$\mathfrak{F}$-subnormal} in
 $G$      if   there is a subgroup
 chain  $$A=A_{0} \leq A_{1} \leq \cdots \leq
A_{n}=G$$  such that  either $A_{i-1} \trianglelefteq A_{i}$ or 
$A_{i}/(A_{i-1})_{A_{i}} \in \mathfrak{F}$ 
  for all $i=1, \ldots , n$.   A formation  $\frak {F}$ is said to be 
$K$-lattice provided in every finite group $G$ the set of all its  
$K$-$\mathfrak{F}$-subnormal subgroups forms a sublattice of the lattice of 
all subgroups of $G$. 

In this paper we consider some new applications of  the theory of  
$K$-lattice  formations.  In particular, we prove the following

{\bf Theorem A.} {\sl Let  $\mathfrak{F}$ be a hereditary $K$-lattice 
  saturated formation containing
 all nilpotent groups. }

(i)  {\sl If every   $\mathfrak{F}$-critical subgroup $H$ of $G$ 
  is $K$-$\mathfrak{F}$-subnormal in $G$ with  
$H/F(H)\in {\mathfrak{F}}$, then 
$G/F(G)\in {\mathfrak{F}}$.}

(ii)  {\sl
 If every  Schmidt  subgroup of $G$     is $K$-$\mathfrak{F}$-subnormal in $G$, then 
   $G/G_{\mathfrak{F}}$ is abelian.}

\end{abstract}


\let\thefootnote\thefootnoteorig

\section{Introduction}

Throughout this paper, all groups are finite and $G$ always denotes
a finite group.  Moreover, $\frak {F}$  is a non-empty class of group, and  
if $1\in \frak {F}$, then 
  $G^{\frak {F}}$  denotes the intersection of all normal
subgroups $N$ of $G$ with  $G/N\in {\frak {F}}$; $G_{\frak {F}}$  is the product
  of all normal
subgroups $N$ of $G$ with  $N\in {\frak {F}}$.

For any   equivalence $\pi$  on the set of all primes $\Bbb{P}$, we write $\text{part} (\pi)$ to 
denote the  partition of $\Bbb{P}$ defined by $\pi$.
 On the other hand, for any partition $\sigma$ of $\Bbb{P}$,
 we write $\text{eq} (\sigma)$ to 
denote the equivalence   on $\Bbb{P}$ defined by $\sigma$.

If $\sigma =\{\sigma_{i} | i\in I\}$ is 
a 
partition of $\Bbb{P}$ (that is,
$\Bbb{P}=\bigcup_{i\in I} \sigma_{i}$ and $\sigma_{i}\cap
\sigma_{j}= \emptyset $ for all $i\ne j$), then $G$ is said to be:   
\emph{$\sigma$-primary} \cite{1} if  $G$ is a $\sigma_{i}$-group for some $i$; 
\emph{$\sigma$-decomposable} (Shemetkov \cite{Shem})  or 
\emph{$\sigma$-nilpotent} (Guo and Skiba  \cite{33}) if $G=G_{1}\times \dots \times G_{n}$ 
for some $\sigma$-primary groups $G_{1}, \ldots, G_{n}$; \emph{$\sigma$-soluble} \cite{1}
 if every chief factor of $G$ is $\sigma$-primary. We use  
${\frak{N}}_{\sigma}$ to denote the class of all $\sigma$-nilpotent  groups.

If  ${\sigma ^{0}}=\{\sigma^{0}_{j}\ |\ j\in J\}$ is another  partition of 
 $\Bbb{P}$, then  we write $ \sigma ^{0} \leq  \sigma $ provided    
$eq (\sigma ^{0}) \subseteq  eq (\sigma )$, that is,  for each $j\in J$ there is $i\in I$
 such that  $\sigma ^{0}_{j}\subseteq \sigma _{i}$.  It is clear that if  
$G$ is $ \sigma ^{0}$-nilpotent (respectively 
\emph{$\sigma ^{0}$-soluble}), then  $G$ is $ \sigma $-nilpotent (respectively 
\emph{$\sigma$-soluble}).

 We say that   $\frak {F}$ is  \emph{$\sigma$-nilpotent} (respectively 
\emph{$\sigma$-soluble}) if every group in  $\frak {F}$ is
\emph{$\sigma$-nilpotent} (respectively \emph{$\sigma$-soluble}). 

 For any  set 
$\{\sigma ^{i} \ |\ i\in I\}  $ 
of 
partitions $\sigma ^{i}$  of $\Bbb{P}$ we put  $$\bigcap_{i\in I}  \sigma 
^{i} = \text{part} (\bigcap_{i\in I}  \text{eq} (\sigma ^{i})).$$

If $\{\sigma ^{i} \ |\ i\in I\}  $ is the set 
of  all 
partitions $\sigma ^{i}$  of $\Bbb{P}$  such that  $\frak {F}$ is
  \emph{$\sigma ^{i}$-nilpotent} (respectively 
\emph{$\sigma ^{i}$-soluble}) for all $i$, then write $\Sigma _{n}(\frak{F})$ (respectively
 $\Sigma _{s}(\frak{F})$) 
to denote  the partition  $\bigcap_{i\in I}  \sigma 
^{i} $.  It is clear that
 $\Sigma _{n}(\frak{F})\in \{\sigma ^{i} \ |\ i\in I\}  $
 (respectively $\Sigma _{s}(\frak{F})\in \{\sigma ^{i} \ |\ i\in I\}  $), 
and   $\Sigma _{n}(\frak{F})$ (respectively $\Sigma _{s}(\frak{F})$) is the smallest element in 
  $\{\sigma ^{i} \ |\ i\in I\}  $, that is, $\Sigma _{n}(\frak{F})\leq \sigma ^{i}$
  (respectively $\Sigma _{s}(\frak{F})\leq \sigma ^{i}$) for all $i$.

Recall that a class  of groups $1\in {\frak {F}}$  is  a \emph{formation}  
if for every group $G$ every homomorphic image of $G/G^{\frak {F}}$
 belongs to $\frak {F}$. The formation $\frak {F}$ is said to be:
\emph{saturated} if $G\in {\frak F}$ whenever $\Phi (G)\leq G^{\frak {F}}$; \emph{hereditary} if
 $H\in \frak {F}$ whenever $H\leq G \in
\frak {F}$.

A  subgroup $A$ of $G$ is said to be  \emph{$\mathfrak{F}$-subnormal in
 $G$ in the sense of Kegel} \cite{55} or  \emph{$K$-$\mathfrak{F}$-subnormal} in
 $G$  \cite[6.1.4]{15}    if   there is a subgroup
 chain  $$A=A_{0} \leq A_{1} \leq \cdots \leq
A_{n}=G$$  such that  either $A_{i-1} \trianglelefteq A_{i}$ or 
$A_{i}/(A_{i-1})_{A_{i}} \in \mathfrak{F}$ 
  for all $i=1, \ldots , n$. In particular, $A$ of $G$ is said to be  
\emph{$\sigma$-subnormal}
in   $G$ \cite{1} provided $A$ is $K$-${\frak{N}}_{\sigma}$-subnormalin $G$, that is, 
   there is a subgroup
 chain  $$A=A_{0} \leq A_{1} \leq \cdots \leq
A_{n}=G$$  such that  either $A_{i-1} \trianglelefteq A_{i}$ or 
$A_{i}/(A_{i-1})_{A_{i}} $  is $\sigma$-primary 
  for all $i=1, \ldots , n$.

  The set  ${\cal L}_{K\frak{F}}(G)$  of all $K$-$\mathfrak{F}$-subnormal
subgroups of  $G$   
is partially ordered with respect to   set inclusion. Moreover, ${\cal L}_{K\frak{F}}(G)$  
   is a lattice  
 since $G\in {\cal L}_{K\frak{F}}(G)  $ and,  by \cite[Lemma 6.1.7]{15}, for any 
 $A_{1}, \ldots , A_{n}\in {\cal L}_{K\frak{F}}(G)$   the subgroup $
 A_{1} \cap  \cdots \cap  A_{n} \in {\cal L}_{K\frak{F}}(G)$, so  this intersection  is the 
  greatest lower   bound for $\{A_{1}, \ldots , A_{n}\}$  in  ${\cal L}_{K\frak{F}}(G)$.

 The formation   $\frak {F}$ is called  \emph{$K$-lattice} \cite{15} if in every 
group $G$ the lattice ${\cal L}_{K\frak{F}}(G)$  is a sublattice of the 
lattice ${\cal L}(G)$
 of all subgroups of $G$. 

 The full classification of  hereditary $K$-lattice saturated formations 
were given in the papers  \cite{22, 333} (see also Ch. 6 in \cite{15}).
The formations of such kind were useful in the study 
of many problems in the theory of finite groups (see, in particular,
 the recent papers
 \cite{Skibaja5, skjpaa, GuoS-JA} and Ch. 6 in \cite{15}). 

In the given paper,  we consider two new  applications of hereditary 
 $K$-lattice  saturated  formations.

 Recall that    $G$ is said to be
 \emph{$\frak {F}$-critical} if $G$ is not
 in $\frak {F}$ but all proper subgroups of $G$
are in $\frak {F}$ \cite[p. 517]{DH}; $G$ is said to be a \emph{Schmodt 
group} provided  $G$ is 
$\frak {N}$-critical, where $\frak {N}$ is the class of all nilpotent groups.

A large number of publications are related to the study of the influence 
on the structure of the group of its critical subgroups,  in particular,  
Schmidt subgroups. It was proved, for example, that if every Schmidt subgroup of 
$G$ is subnormal, then $G'\leq F(G)$  \cite{Sem, KM}. Later, this result 
was generalized in the paper \cite{Kh1}, where it was proved that
if every Schmidt subgroup of 
$G$ is $\sigma$-subnormal in $G$, then $G'\leq F_{\sigma}(G)$ (here  $F_{\sigma}(G)=G_{{\frak{N}}_{\sigma}}$ is the
 \emph{$\sigma$-Fitting subgroup} of $G$, that is,
  the product of all normal $\sigma$-nilpotent subgroups
 of $G$).

 Our first observation is the following generalization of these results.  

{\bf Theorem A.} {\sl Let  $\mathfrak{F}$ be a hereditary $K$-lattice 
  saturated formation containing
 all nilpotent groups. }

(i)  {\sl If every   $\mathfrak{F}$-critical subgroup $H$ of $G$ 
  is $K$-$\mathfrak{F}$-subnormal in $G$ with  
$H/F(H)\in {\mathfrak{F}}$, then 
$G/F(G)\in {\mathfrak{F}}$.}

(ii)  {\sl
 If every  Schmidt  subgroup of $G$     is $K$-$\mathfrak{F}$-subnormal in $G$, then 
  $G/G_{\mathfrak{F}}$ is abelian.}

 Note that if $\mathfrak{F}=\mathfrak{N}$  is the formation of all 
nilpotent groups, then  a subgroup $A$ of $G$ $K$-$\mathfrak{F}$-subnormal in 
$G$  if and only if $A$ is subnormal in $G$. Hence we get from Theorem 
A(i) the following two known results.

{\bf Corollary 1.1} (Semenchuk \cite{Sem}).  {\sl  If 
  every Schmidt subgroup of $G$ is
 subnormal in $G$, then $G$ is metanilpotent.}

{\bf Corollary 1.2} (Monakhov and Knyagina \cite{KM}).
  {\sl   If                        
  every Schmidt subgroup of $G$ is
 subnormal in $G$, then $G/F(G)$ is abelian.}

From Theorem A(ii) we get the following

{\bf Corollary 1.3 } (Al-Sharo,  Skiba \cite{Kh1}). {\sl If every Schmidt subgroup of $G$ is
 ${\sigma}$-subnormal in $G$, then $G/F_{\sigma}(G)$ is  abelian}.

  Recall that if $M_n < M_{n-1} < \ldots < M_1 < M_{0}=G$ (*), where
 $M_i$ is a maximal subgroup of  $M_{i-1}$ for all $i=1, \ldots ,n$,
 then the chain   (*)  is said to be  a \emph{maximal chain of $G$ of
 length $n$} and $M_n $ ($n > 0$),  is an \emph{$n$-maximal subgroup}  of $G$.

 If $\mathfrak{F}$ is a saturated formation containing all nilpotent groups,
 then $G\in \frak{F}$ if and only
 if every maximal subgroup of $G$ is $K$-$\mathfrak{F}$-subnormal in 
$G$.  But  when we deal with  hereditary $K$-lattice  saturated formations, the following
 result is true.

{\bf Theorem  B}.  {\sl  Let  $\mathfrak{F}$ be a hereditary $K$-lattice 
  saturated formation containing
 all nilpotent groups and $\sigma =\Sigma _{s}(\frak{F})$. Then the following statements hold:
  }

(i)  {\sl Every maximal
 chain of $G$ of  length $2$  includes a proper 
$K$-$\mathfrak{F}$-subnormal subgroup of $G$  if and only if either $G\in \frak{F}$ 
or $G\not \in \frak{F}$ is a Schmidt group  with abelian Sylow subgroups.}

(ii) {\sl If  every maximal
 chain of $G$ of  length $3$ includes a proper 
$K$-$\mathfrak{F}$-subnormal subgroup of $G$, then  $G$ is $\sigma$-soluble.  }

In the case $\mathfrak{F}=\mathfrak{N}$  we get from Theorem 
B the following two known results.

{\bf Corollary 1.4} (Spencer \cite{spen}). {\sl If every
maximal chain  of $G$  of length 3    includes a proper
 subnormal subgroup of $G,$ then $G$ is soluble.   }

{\bf Corollary 1.5} (Spencer \cite{spen}). {\sl If every
maximal chain  of a non-nilpotent group $G$  of length 2    includes a proper
 subnormal subgroup of $G,$ then $G$ is  a Schmidt group  with abelian Sylow subgroups.   }
 
In the next two results, $\pi$ is a non-empty set of primes.

{\bf Corollary 1.6.} {\sl Suppose that $\frak{F}=\frak{F}_{\pi'}$ is
 the class of all $\pi'$-groups.
 If every
maximal chain  of $G$  of length 3    includes a proper
$K$-$\frak{F}$-subnormal subgroup of $G,$ then $G$ is $\pi$-soluble.   }

{\bf Corollary 1.7.} {\sl Suppose that $\mathfrak{F}=\mathfrak{N}_{\sigma}$, where
 $\sigma =\{\pi, \pi'\}$.
 If every
maximal chain  of $G$  of length 3    includes a proper
$K$-$\frak{F}$-subnormal subgroup of $G,$ then $G$ is $\pi$-separable.   }

\section{Proof of Theorem A}

In view of Proposition 2.2.8  in \cite{15}, we have the following

{\bf Lemma 2.1.}    {\sl Let    
$\mathfrak{F}$ be  a non-empty formation. If  $N$ and $U$ are  subgroups of 
$G$ such that $N$ is  normal in $G$  and $G=NU$, then:   }

(i)  {\sl $(G/N)^{\frak{F}}=G^{\frak{N}}N/N,$ and} 

 (ii) {\sl $G^{\frak{F}}N=U^{\frak{F}}N$.}

{\bf Lemma 2.2} (See Corollary 4.2.1 in \cite{Shem}). {\sl If   $\mathfrak{F}$ is a 
 saturated   formation containing
 all nilpotent groups  and    $E$ is a normal subgroup of $G$ such that $E/E\cap\Phi (G)
 \in \mathfrak{F}$, then $E\in \mathfrak{F}$.}

Let $\mathfrak{F}$     is a hereditary formation. Then
 a  subgroup $A$ of $G$ is said to be  \emph{$\mathfrak{F}$-subnormal in
 $G$}  if   there is a subgroup
 chain  $$A=A_{0} \leq A_{1} \leq \cdots \leq
A_{n}=G$$  such that   
$A_{i}/(A_{i-1})_{A_{i}} \in \mathfrak{F}$ 
  for all $i=1, \ldots , n$.   A formation  $\frak {F}$ is said to be 
\emph{lattice} \cite{15} provided in every  group $G$ the set of all its  
$\mathfrak{F}$-subnormal subgroups forms a sublattice of the lattice of 
all subgroups of $G$.

{\bf Lemma 2.3.} {\sl If   $\mathfrak{F}$ is a hereditary
 $K$-lattice  saturated formation containing
 all nilpotent groups and    $G$ is  an  $\mathfrak{F}$-critical group with  
$G/F(G)\in {\mathfrak{F}}$, then $G$ is a Schmidt group.}

{\bf Proof.  }  Suppose that this lemma is false. By 
 \cite[Theorem  6.3.15]{15}, $\frak{F}$ is  lattice.  
Hence by Lemma 13 in \cite{333},  $(G/\Phi(G))^{\mathfrak{F}}$ is a  
non-abelian minimal normal subgroup of  $G/\Phi(G)$. But Lemma 2.1 implies 
that $$(G/\Phi(G))^{\mathfrak{F}}=G^{\mathfrak{F}}\Phi(G)/\Phi(G)\simeq 
G^{\mathfrak{F}}/(G^{\mathfrak{F}}\cap \Phi(G))$$  is nilpotent. Hence 
$(G/\Phi(G))^{\mathfrak{F}}$ is nilpotent by Lemma 2.2, so it is abelian. This contradiction 
completes the proof of the lemma.

{\bf Lemma 2.4} (See \cite[Theorems  6.3.3, 6.3.15]{15}).
 {\sl Let  $\mathfrak{F}$ be a hereditary $K$-lattice 
  saturated formation containing
 all nilpotent groups and $\sigma =\{\sigma _{i}\ |\ i\in I\}=\Sigma _{n}(\frak{F})$.
 Then the following statements hold:  }

(i) {\sl If $\pi \subseteq\sigma _{i}$ for some $i$, then each soluble $\pi$-group 
 is contained in $\mathfrak{F}$. }

(ii) {\sl If $A\in \frak{F}$ and $B\in \frak{F}$ are 
 $K$-$\frak{F}$-subnormal subgroups of $G$, then 
$\langle A, B \rangle \in \frak{F}$.}

{\bf Lemma 2.5} (See \cite[Theorem  A]{15}).   {\sl  If $G$ is $\sigma$-soluble,
 for some partition 
 $\sigma =\{\sigma_{i} | i\in I\}$  of $\Bbb{P}$, then
 $G$ possesses a Hall $\sigma _{i}$-subgroup for all $i$.}

We say that $G$ is \emph{$\sigma$-metanilpotent} if $G/ F_{\sigma}(G)$ is $\sigma$-nilpotent.

{\bf Lemma 2.6} (See \cite[Proposition 4.2]{Kh1}). {\sl If $A$ is a
 $\sigma$-subnormal subgroup of $G$, for some partition  $\sigma$ of $\Bbb{P}$, and $A$ is
 $\sigma$-metanilpotent (respectively $\sigma$-nilpotent),  then $A^{G}$ is
 $\sigma$-metanilpotent (respectively $\sigma$-nilpotent). }
 
 \

 {\bf Lemma 2.7} (see \cite[Lemma 2.6]{1}). {\sl If $A$ is a
 $\sigma$-subnormal subgroup of $G$, where $\sigma =\{\sigma_{i} | i\in I\}$,  and $A$ is
 $\sigma _{i}$-group,  then $A\leq O_{\sigma _{i}}(G)$. }

 {\bf Proof of Theorem A.}   Suppose
 that this theorem is false and let $G$ be a counterexample
 of minimal order.  Then $G\not \in \mathfrak{F}$.

(*) {\sl If $H$ is an $\mathfrak{F}$-critical subgroup of $G$, then
 $H^{\mathfrak{F}}\leq F(G)$. Hence $G$ possesses  an abelian minimal normal
 subgroup,   $R$ say.}

Since $G\not \in \mathfrak{F}$, $G$ has an 
 $\mathfrak{F}$-critical subgroup, $H$ say. The hypothesis implies that 
$H/F(H)\in {\mathfrak{F}}$ and  so  $H^{\mathfrak{F}}\leq F(H)$ since 
$\mathfrak{F}$.  Moreover, 
$H^{\mathfrak{F}}$ is subnormal in $G$ by \cite[Lemma 6.1.9]{15}.
 But then, by using
 \cite[Theorem  
6.3.3]{15} in the case when $\frak{F}=\frak{N}$, we get that   
$1  < H^{\mathfrak{F}}\leq F(G)$.  Hence we have (*).

(i)       Assume that this  is false.

(1)  {\sl The hypothesis  holds for every subgroup of $G$.
 Hence $E/F(E)\in {\mathfrak{F}}$ for every proper
 subgroup $E$ of $G$. }

If $E\in  \frak{F}$, it is clear. Now assume that   $E\not \in {\mathfrak{F}}$ and let  
$H$ be any $\mathfrak{F}$-critical subgroup  of $E$, then $H$ is   
$K$-$\mathfrak{F}$-subnormal in $G$  by hypothesis, so  $H$ is  
$K$-$\mathfrak{F}$-subnormal in $E$   by \cite[Lemma 6.1.7]{15}. Therefore the 
hypothesis holds for $E$, so the choice of $G$ implies that  we have  (1).

(2) {\sl $(G/N)/F(G/N)\in {\mathfrak{F}}$ for  each  minimal normal subgroup
 $N$ of $G$.   }

In view of the choice of $G$, it is enough to show that the hypothesis 
holds for $G/N$.  Suppose that this is false. Then $G/R\not \in 
\mathfrak{F}$. 
Let  $K/N$ be any $\mathfrak{F}$-critical subgroup  of $G/N$, and let $L$ be
 any minimal supplement to $N$ in $K$.  Then $L\cap N\leq \Phi (L)$. 
Moreover, $K/N=LN/N\simeq L/L\cap N$ is an $\mathfrak{F}$-critical group, so $L/\Phi (L)$ 
is an $\mathfrak{F}$-critical group.  Now let  $A_{1}, 
\ldots , A_{t}$  be the set of all  $\mathfrak{F}$-critical subgroups of 
$L$ and  $V= \langle A_{1}, 
\ldots , A_{t}\rangle $.  First we show that $L=V$. It is clear that  $V$ is 
normal in $L$. Suppose that $V < L$. Then $V\Phi (L) < L$, so  $V\Phi (L)/\Phi (L) < L/\Phi 
(L)$ and hence    $$V\Phi (L)/\Phi (L)\simeq V/V\cap \Phi (L)\in {\frak{F}}$$ 
 since $L/\Phi (L)$ 
is an $\mathfrak{F}$-critical group.   
  But then $V\in \frak{F}$ by Lemma 2.2,
 so  $A_{1}\in \frak{F}$  since $\frak{F}$ is hereditary by hypothesis. 
This contradiction shows that   $V=L=\langle A_{1}, 
\ldots , A_{t}\rangle$. Since $\frak{F}$ 
is a $K$-lattice formation, it follows that $L$ is  
$K$-$\mathfrak{F}$-subnormal in $G$. Hence $K/N=LN/N$ is 
$K$-$\mathfrak{F}$-subnormal in $G/N$ by  \cite[Lemma 6.1.6]{15}.

Finally, we show that  $(K/N)/F(K/N)\in {\mathfrak{F}}$.  
Claim (*) implies that  for each  $i$ we have  $A_{i}^{\frak{F}} \leq F(G)\cap L\leq F(L)$.
 Hence  $$A_{i}F(L)/F(L)\simeq 
A_{i}/A_{i}\cap F(L)\simeq (A_{i}/A_{i}^{\frak{F}})/((A_{i}\cap F(L))/A_{i}^{\frak{F}})\in
 {\frak{F}}.$$  
On the other hand,   $A_{i}F(L)/F(L)$ is   
$K$-$\mathfrak{F}$-subnormal in $L/F(L)$ by \cite[Lemma 6.1.6]{15}.   Hence $$L/F(L)= 
\langle A_{1}F(L)/F(L), 
\ldots A_{t}F(L)/F(L)\rangle \in 
\frak{F}$$ by Lemma 2.4(ii).  
 Thus $L^{\frak{F}}\leq F(L)$. 
Therefore, by Lemma 2.1(ii),    
   N $$(K/N)^{\frak{F}}=K^{\frak{F}}N/N=L^{\frak{F}}N/N\simeq L^{\frak{F}}/L^{\frak{F}}\cap
 N$$ is nilpotent.  Therefore  $(K/N)/F(K/N)\in {\mathfrak{F}}$, so the 
hypothesis holds for $G/N$.

(3) {\sl $R$ is the  unique  minimal normal subgroup of $G$ and 
 $R\nleq \Phi (G)$. Hence $R=F(G)=O_{p}(G)=C_{G}(R)$ for some prime $p$.}

In view of Claim (2)   and Lemma 2.1, we get that  $$(G/R)^{\mathfrak{F}}=DR/R\simeq 
D/D\cap R$$ is nilpotent.     Hence the choice of $G$ and Lemma 2.2 imply  that 
$R\nleq \Phi (G)$. Finally, if $G$ has a  minimal normal subgroup $N\ne R$,
 then $D\simeq
 D/1=D/R\cap N$ is nilpotent and so $G/F(G) \in \mathfrak{F}$, contrary 
 the choice of $G$.  Hence $R$ is the unique  minimal normal subgroup of $G$, so
 we have (3) by \cite[A, 15.6]{DH}.

(4) {\sl  Final contradiction for (i).} 

Assume that $G/R\not \in {\mathfrak{F}}$. Then, in view of Claim (3), for some maximal subgroup 
$M$ of $G$ we have $G=R\rtimes M$, where $M\simeq G/R\not \in \frak{F}$. 
Now let $H$ be   any $\mathfrak{F}$-critical subgroup  of $M$. Then 
$1  < H^{\frak{F}}\leq F(G)\cap M$ by Claim (*).  But  $F(G)\cap M=R\cap M=1$  by Claim (3).
 Therefore $H^{\frak{F}}=1 $ and so $H\in \mathfrak{F}$  , a contradiction.  Thus Statement  (i) is true.

(ii)   Assume that this assertion is false.

(5) {\sl The hypothesis holds for every subgroup of $G$.
 Hence $E/E_{\frak{F}}$
 is abelian for every proper
 subgroup $E$ of $G$  } (See Claim (1)).

(6) {\sl $(G/N)/(G/N)_{\frak{F}}$ is abelian
 for  each  minimal normal subgroup  $N$ of $G$}  (See Claim (2)).

(7) {\sl $R$ is the unique  minimal normal subgroup of $G$ and 
 $R\nleq \Phi (G)$. Hence $G=R\rtimes M$ and $R=F(G)=O_{p}(G)=C_{G}(R)$ for some prime $p$
 and 
  some maximal subgroup $M$ of $G$.}

First note that $G'=G^{\frak {A}}$, where   $\frak {A}$ is the class of 
all abelian groups.  
Claim (6)   and Lemma 2.1 imply that $$(G/R)'=G'R/R\simeq 
G'/G'\cap R \in \frak{F}. $$     Hence the choice of $G$ and Lemma 2.2 imply  that 
$R\nleq \Phi (G)$. Finally, if $G$ has a  minimal normal subgroup $N\ne R$,
 then  $G'\simeq
 G'/1=G'/R\cap N \in \frak{F}$ and so $G/G_{\frak{F}}$ is abelian,
 contrary to the choice of $G$. Hence we have (7) by \cite[A, 15.6]{DH}.

(8) {\sl $G/R \in {\mathfrak{F}}$.}   (See the proof of Claim (4)   and use Claim (7)).

(9) {\sl   If $\sigma =\{\sigma_{i} | i\in I\}$, where $\sigma=\Sigma _{n} 
(\frak{F})$, then  $G$ is $\sigma$-soluble and 
  $$ R= O_{\sigma _{i}}(G)=F_{\sigma}(G)$$ is a Hall $\sigma  _{i}$-group 
of $G$ for some $i$. }

Since $G/R\in {\mathfrak{F}}$ by Claim (8) and also    
every group in ${\mathfrak{F}}$ is $\sigma$-nilpotent by definition (and so 
$\sigma$-soluble),   Claim (7) implies that $G$ is $\sigma$-soluble.

Claims (7),  (8) and Lemma 2.5 imply that    $R\leq F_{\sigma}(G)=O_{\sigma _{i}}(G) \leq H$,
 where $H$    is a Hall $\sigma  _{i}$-group of $G$ for some $i$.     
 Since  $R=G^{\frak{F}}$ by Claim (8) and every group in $\frak{F}$ is
 $\sigma$-nilpotent,   
 $H/O_{\sigma _{i}}(G)$ is a normal subgroup of 
$G/O_{\sigma _{i}}(G)$. Hence  $H\leq O_{\sigma _{i}}(G)$, so 
$ H=O_{\sigma _{i}}(G)$. Therefore  
$F_{\sigma}(G)=H_{i}$. Now assume that $R <  H$. 
 By the   Schur-Zassenhaus  theorem,  
 $G$ has a  ${\sigma _{i}}$-complement, $W$ say. Then 
$V=RW < G$, so Claim (5) implies that  $V/V_{\frak{F}}$ is abelian.
 Since $R$ is a Hall $\sigma  _{i}$-group of $V$, 
 $V_{\frak{F}}=R\times (V_{\frak{F}}\cap W)$. 
But Claim (7) implies that $C_{G}(R)\leq R$. Hence $V_{\frak{F}}=R$ and so 
$W\simeq V/R$ is abelian.

Suppose that  $H\not \in \frak{F}$  and let $A$ be an $\frak{F}$-critical 
subgroup of $H$. Then,   by Lemma 2.3 and \cite[III, 5.2]{hupp},  $A=A_{q}\rtimes A_{r}$ is a Schmidt group,
 where $A_{q}\in \text{Sylow}_{q}(A)$, 
  $A_{r}\in \text{Sylow}_{r}(A)$ and  $q, r \in \sigma _{i}$. Hence  
 $A\in \frak{F}$ by Lemma 2.4(i). This contradiction shows that  $H \in \frak{F}$, so  $H 
\leq  G_{\frak{F}}$.  
 Therefore $G/G_{\frak{F}}=(G/H)/(G_{\frak{F}}/H)$
 is abelian since $G/H\simeq W$ is abelian, contrary to the choice  of $G$. 
Thus $R=H$, so we  have (9).

(10) {\sl $M$   is a Miller-Moreno group (that is, a $\frak{U}$-critical
  group, where $\frak{U}$ is
 the class of all abelian groups).  Moreover, $M$ is a $q$-group for some prime $q\ne p$.} 

First note that $M$ is a  Hall  ${\sigma _{i}}'$-subgroup $M$ of $G$  by 
Claims (7) and (9). 
Now, let $S$ be any maximal subgroup of $M$. Then $RS/(RS)_{\frak{F}}$ is 
abelian    by Claim (5). In view of Claim (7), 
$R= (RS)_{\frak{F}}$ and hence $S\simeq RS/R$ is abelian.   
 Therefore the    
 choice of $G$ and Claim (7) imply that  $M$ is a $\frak{U}$-critical
  group.   Therefore,  $M$ is either a  Schmidt  group 
 or a minimal non-abelian group of prime power order $q^{a}$. 
In the former case, by \cite[III, 5.2]{hupp},   $M=Q\rtimes V$,
 where $Q=M^{\frak{N}}\in \text{Sylow}_{q}(M)$, $V\in \text{Sylow}_{r}(M)$ and $q\ne r$. 
Since  $ R=C_{G}(R)$ is a Hall $\sigma  _{i}$-group of $G$  by Claim (9), 
$RQ\not \in \frak{F}$, so $RQ$ has an $\frak{F}$-critical subgroup $A$. 
       By 
hypothesis, $A$ is $K$-$\frak{F}$-subnormal in $G$, so it is 
$\sigma$-subnormal in $G$. Therefore $A^{G}$ is a meta-$\sigma$-nilpotent 
group by Lemma 2.6.  
 Similarly,  $RV$ has an $\frak{F}$-critical subgroup $B$ and
  $B^{G}$  is  
meta-$\sigma$-nilpotent. But then  $G=A^{G}B^{G}$ is a 
meta-$\sigma$-nilpotent. Indeed, $F_{\sigma}(A)$ is a characteristic
 $\sigma$-nilpotent subgroup of $A$, so $F_{\sigma}(A)$  is a 
 $\sigma$-nilpotent $\sigma$-subnormal subgroup of $G$. Hence
 $F_{\sigma}(A)\leq R=F_{\sigma}(G)$ and $F_{\sigma}(G)\cap A\leq 
F_{\sigma}(A)$, which implies that 
 $AF_{\sigma}(G)/F_{\sigma}(G)\simeq A/A\cap F_{\sigma}(G)$ is $\sigma$-nilpotent $\sigma$-subnormal 
subgroup of $G/F_{\sigma}(G)$. Similarly,  $BF_{\sigma}(G)/F_{\sigma}(G)$ is a 
 $\sigma$-nilpotent $\sigma$-subnormal 
subgroup of $G/F_{\sigma}(G)$. Therefore $G/R=G/F_{\sigma}(G)\simeq M$ is  
$\sigma$-nilpotent by Lemma 2.6. Therefore, since    $M$ is $K$-$\frak{F}$-subnormal in $G$ 
by hypothesis,  $M$ is $\sigma$-subnormal in $G$ and  so $M\leq F_{\sigma}(G)=R$ by
 Lemma 2.6. This contradiction shows that  we have the second case and so Claim (10) holds.

{\sl Final contradiction for (ii). }  From Claims  (7), (9) and (10) we get that 
$G=R\rtimes M$, where $R=G^{\frak{F}}$ and $M$ is a $\frak{U}$-critical
   $q$-group for some prime $p\in \sigma _{i}$ and $q\in \sigma _{j}$, 
where $i\ne j$. Let $Z$ be a group of order $q$ in $Z(M) $  and  $E=RZ$.
Then $E\not \in \frak{F}$ by Claim (7) since  all $\{p, q\}$-groups  contained
 in  $ \frak{F}$ are nilpotent. 
Let   $A=A_{p}\rtimes Z$ be an $ \frak{F}$-critical subgroup of $E$.

Note  that $R=R_{1}\times \cdots    \times R_{t}$, where  $R_{k}$ is 
a minimal normal subgroup of $E$ for all $k=1, \ldots, t$   by the Mashcke's
 theorem.  
 Suppose that  $A < E$. Then  there is a  proper
 subgroup $V$ of $E$ such that $A\leq V$ 
and either $E/V_{E}$ is a $p$-group or $V$ is normal in $E$.
Then   $Z\leq  V_{E} < E$, so for some $k$ we have $R_{k}\nleq  V_{E}$.    Hence 
$R_{k}\leq C_{E}(V_{E})$, so $R_{k} \leq N_{G}(Z)=M$. Thus $Z$ is normal in $G$ since $M$ is a 
maximal subgroup of $G$ and hence $Z\leq C_{G}(R)=R $,   
 a contradiction. Therefore $E=A$, so $R=P$ and $Z$ acts irreducibly on 
$R$.

  It is clear that  $Z\leq  \Phi (M)$ and  so  
 every maximal subgroup  of $W$ acts irreducibly on $R$, which implies that every
 maximal subgroup of 
$W$ is cyclic. Hence $q=2$ and so $|R|=p$. It follows that $G/R=G/C_{G}(R)\simeq W$ is abelian, 
contrary to Claim (10). Thus Statement  (ii) is true. 

The theorem is proved.

\section{Proof of Theorem  B}

    If  $G \in \frak{F}$, then every subgroup of $G$ is clearly 
 $K$-$\mathfrak{F}$-subnormal in $G$  since  $\mathfrak{F}$ is  a 
hereditary  formation by hypothesis. Moreover, if $G$  is a Schmidt group
  with abelian Sylow subgroups,   then  every proper subgroup $H$ of $G$ is 
subnormal in $G$ and  so  it is $K$-$\mathfrak{F}$-subnormal in $G$.

Now assume that $G\not \in \mathfrak{F}$ and that   every maximal
 chain of $G$ of  length $2$  includes a proper 
$K$-$\mathfrak{F}$-subnormal subgroup of $G$. We show that in this case  
$G$ is a Schmidt group  with abelian Sylow subgroups.
First note that   $G\not \in \frak{F}$ implies that  for some maximal subgroup $M$ 
of $G$ we have   $G/M_{G}\not \in \mathfrak{F}$ since $\mathfrak{F}$ is a 
  saturated formation.  Hence $M$ is not $K$-$\mathfrak{F}$-subnormal in $G$.
 Therefore  every  maximal subgroup of $M$ is  
$K$-$\mathfrak{F}$-subnormal in $G$ by hypothesis.  Hence  
 $M$ is a  cyclic Sylow $p$-subgroup of $G$ because $\mathfrak{F}$ is a 
 lattice  formation.  Therefore $G$ is soluble by
 the Deskins-Janko-Thompson theorem \cite[7.4, IV]{hupp}.
 Suppose that $M_{G}\ne 1$ and let $R$ be a minimal normal subgroup of $G$ 
contained in $M_{G}.$ Then $R\leq Z(G)$, since $M\leq C_{G}(R)$ and $M$ is 
a maximal subgroup of $G$ and it  is evidently  not normal in $G$. In view of 
\cite[Lemma 6.1.6]{15},  the hypothesis holds for $G/R$ and, clearly,
 $G/R\not \in \mathfrak{F}$. Hence  $G/R$ is  a Schmidt
 group with abelian Sylow subgroups.  Therefore, if $V/R$ is any 
maximal subgroup of $G/R$, then $V/R$  is nilpotent and so $V$ is 
nilpotent since $R\leq Z(G)$. Therefore $G$ is a Schmidt group. It is clear also that  
the  Sylow subgroups of $G$ abelian. 

Now assume that $M_{G}=V_{G}=1$. Then $G=R\rtimes M$,
 where $R$ is a
 minimal normal subgroup of $G$ and $R$ is a Sylow $p$-subgroup of $G$ for some prime $q\ne p$.
 Let $V$ be the maximal subgroup of $M$.  Then $V, R, M\in \frak{F}$ since $\frak{F}$ 
contains all nilpotent group by hypothesis and  so $G$ is an $\frak{F}$-critical group because 
$RV\in \frak{F}$ by Lemma 2.4(ii). But then $G$ is a Schmidt group  by 
Lemma 2.3, and the   
Sylow subgroups  $R$ and $M$ of $G$ are abelian.
  Thus Statement (i) is true.

(ii)    Suppose that this assertion  is false and let $G$ be a counterexample of minimal
 order.  Then $G\not \in \frak{F}$, since otherwise $G$ is 
$\sigma$-soluble by  definition of $\sigma$.  
 Let   $\Sigma 
_{n}(\frak{F})=\sigma ^{0}=\{\sigma_{j}^{0} | j\in J\}$.
 Then, evidently, $\sigma \leq \sigma ^{0}$.

  First note   that  $G/R$  is
$\sigma$-soluble. Indeed, if $R$ is a maximal subgroup or a 2-maximal
subgroup of $G$, it is clear. Otherwise,   the hypothesis
holds for $G/R$ by \cite[Lemma 6.1.6]{15}, so the choice of $G$ implies that  $G/R$ is
$\sigma$-soluble. Hence $R$ is the unique minimal normal subgroup of $G$   and $R$ is not $\sigma$-soluble. Hence
 $R$ is not abelian and $R\leq G^{\frak{F}}$.

 Let $p$ be any odd prime dividing $|R|$ and $R_{p}$ a Sylow $p$-subgroup of $R$.
 The Frattini argument  implies  that there is a maximal subgroup  $M$
 of $G$ such that $N_{G}(R_{p})\leq M$ and $G=RM$. It is clear that $M_{G}=1$, so $M$
  is not $K$-$\frak{F}$--subnormal in $G$ since $G/M_{G}\simeq G$.
   Let $D=M\cap R$. Then $R_{p}$ is a Sylow $p$-subgroup of $D$.

(1) {\sl $D$ is not nilpotent.   Hence $D\nleq  \Phi (M)$ and $D$ is not a $p$-group.}

Assume that  $D$ is a nilpotent. Then
 $R_{p}$ is normal in $M$. Hence $Z(J(R_{p}))$ is normal
 in $M$. Since $M_{G}=1$, it follows that $N_{G}(Z(J(R_{p})))=M$ and so
 $N_{R}(Z(J(R_{p})))=D$   is  nilpotent. This implies that $R$ is $p$-nilpotent
 by Glauberman-Thompson's theorem on
 the normal $p$-complements. But
then $R$ is a $p$-group, a contradiction. Hence we have (1).

(2) $R < G$.

Suppose that $R=G$ is a simple non-abelian group. Assume that some proper
 non-identity subgroup $A$ of $G$ is $K$-$\frak{F}$-subnormal in $G$. Then there
 is a subgroup chain
 $A=A_{0} \leq A_{1} \leq \cdots \leq A_{n}=G$  such that  either $A_{i-1}\trianglelefteq
  A_{i}$ or
$A_{i}/(A_{i-1})_{A_{i}}\in \frak{F}$ for all $i=1, \ldots ,t$.
  Without loss of generality, we can assume that $M=A_{n-1} < G$. Then $M_{G}=1$
 since  $G=R$ is simple, so $G\simeq G/1\in \frak{F}$,   
 a contradiction. Hence every proper $K$-$\frak{F}$-subnormal
 subgroup of $G$ is trivial.

Let $Q$ be a Sylow $q$-subgroup of $G$,  where $q$ is the smallest prime dividing $|G|$,
  and let $L$ be a maximal subgroup  of $G$ containing $Q$.
Then, in view of \cite[IV, 2.8]{hupp}, $|Q| > q$.  Let $V$ be a maximal subgroup of $Q$.
 If $|V|=q$, then $Q$ is abelian, so $Q < L$ by \cite[IV, 7.4]{hupp}. Hence there
 is a 3-maximal subgroup $W$ of
 $G$ such that $V\leq W$. But then some proper non-identity  subgroup of $G$ is
 $K$-$\frak{F}$-subnormal in $G$ by hypothesis, a contradiction. Therefore $|V| > q$,
 which again
 implies that some proper non-identity  subgroup of $G$ is $K$-$\frak{F}$-subnormal
 in $G$.  This contradiction shows that we have (2).

(3) {\sl $M$ is $\sigma$-soluble.}

If every maximal subgroup of $M$ has prime order, it is evident. 
Otherwise, let $L < T < M$, where $T$ is a maximal subgroup of $M$ and $L$ is a
 maximal subgroup of $T$.
Since $M$ is not $K$-$\frak{F}$-subnormal in $G$, either $L$ or $T$ is
 $K$-$\frak{F}$-subnormal in $G$ and so it is $K$-$\frak{F}$-subnormal in $M$  by
 \cite[Lemma 6.1.7]{15}. Hence  $M$ is $\sigma$-soluble by Part (i).

(4) {\sl $M=D\rtimes T$, where $T$ is a maximal   subgroup of $M$ of  prime order.}

In view of Claim (1), there is a maximal subgroup $T$ of $M$ such that $M=DT$.
 Then $G=RM=R(DT)=RT$ and so, in view of Claim (2), $T\ne 1$. Hence $G$ 
has no a proper subgroup $V$ such that either  $V \trianglelefteq G$ or 
$V/V_{G} \in \mathfrak{F}$.  
Therefore $T$ is not  $K$-$\frak{F}$-subnormal in $G$.

 Assume that $|T|$ is not a prime and let $V$ be a
 maximal subgroup of $T$. Since $M$ and $T$ are  not 
$K$-$\frak{F}$-subnormal in $G$, every maximal subgroup  of $T$ is $K$-$\frak{F}$-subnormal
 in $G$ and so  it is also 
 $K$-$\frak{F}$-subnormal in $T$.  Then  $T, V \in \frak{F}$.
  Moreover, since $|T|$ is not a prime and  $V\ne 1$. 
 Claim (3) implies that  
    $V$ is  $\sigma$-soluble.      Then $V$ is $\sigma ^{0}$-soluble,
 so  
 for some $i$ we have $O_{\sigma _{i}^{0}}(V)\ne 1$. On the other hand,  
$O_{\sigma _{i}^{0}}(V)\in \mathfrak{F}$ since $V\in 
\frak{F}$ and  $O_{\sigma _{i}^{0}}(V)$  
 $K$-$\frak{F}$-subnormal in $G$ since  $V$ is   
$K$-$\frak{F}$-subnormal in $G$. Therefore    $R\leq (O_{\sigma 
_{i}^{0}}(V))^{G} \in \frak{F}$  by Lemma 2.4(ii). Hence $R$ is 
$\sigma$-soluble since every group in $ \frak{F}$ is $\sigma$-soluble by definition.   This 
contradiction completes the proof of  (4).

{\sl Final  contradiction for (ii). } Since $T$ is a maximal subgroup of $M$ and
 it is cyclic, $M$ is soluble and so $|D|$ is a prime power, which contradicts (1).
 Thus Statement  (ii) is true.

 The theorem is proved.

\end{document}